\documentclass[11pt]{article}

\usepackage{graphicx}
\usepackage{url}
\usepackage{natbib}
\usepackage{amssymb}
\usepackage{mathptmx}  % use Times fonts if available on your TeX system
\usepackage[T1]{fontenc}

\newcommand {\myfigure} [1] {\textbf {Figure #1}}
\newcommand {\myfigures} [1] {\textbf {Figures #1}}
\newcommand {\mytable} [1] {\textbf {Table #1}}
\newcommand {\mytables} [1] {\textbf {Tables #1}}
\newcommand {\first} {\vrule depth0ex height2.125ex width0pt}

\newcommand {\chk} {\checkmark}

\begin{document}

\title{Mathematics Turned Inside Out:\\ The Intensive Faculty Versus the Extensive Faculty\footnote {Please cite the article that will appear in \textsl {Higher Education: The International Journal of Higher Education and Educational Planning\/} and which can be accessed online through the DOI link on the arXiv page for this document.}}

\author{Joseph F. Grcar\,\thanks{6059 Castlebrook Drive; Castro Valley, CA 94552-1645 USA} \,\footnote {jfgrcar@comcast.net, jfgrcar@gmail.com.}}

\date {}

\maketitle

\begin{abstract}
Research universities in the United States have larger mathematics faculties outside their mathematics departments than inside. Members of this ``extensive'' faculty conduct most mathematics research, their interests are the most heavily published areas of mathematics, and they teach this mathematics in upper division courses independent of mathematics departments. The existence of this de facto faculty challenges the pertinence of institutional and national policies for higher education in mathematics, and of philosophical and sociological studies of mathematics that are limited to mathematics departments alone.
\end{abstract}

\bigskip

\noindent
\textbf {Keywords}: bibliometrics, mathematics, research university, sociology of science


\section {Introduction}

This paper provides a quantitative basis for exploring policy issues about mathematics education, by painting a comprehensive picture of how, what, and where mathematics is currently studied in American research universities. The data show that mathematics departments primarily teach ``service'' courses for general education while conducting research on subjects unrelated to undergraduate education. Other departments teach and publish papers on the mathematics subjects that upper division undergraduates study most heavily. 

Mathematics has two qualities that make it unique in universities. First, the eponymous department usually has the most student-contact hours of any. This instruction is concentrated in service courses for liberal education requirements or for prerequisites in other fields. Both groups of students place heavy burdens on the service courses: the first by their quantity, the second because the quality of instruction may determine success in the other fields. However, the mathematics faculty who supervise the service courses do not perform quantitative analyses for the most part; indeed, their research topics need have no physical realization (philosophers discuss the epistemology of such knowledge as the question of mathematics ``foundations''). Second, faculty outside mathematics departments provide most of the specialized mathematics eduction beyond the service courses. They also write the majority of mathematics papers, moreover, on subjects relevant to advanced undergraduate instruction.  Since mathematics is not the focus of a single academic department, research universities thus are witness to what is best described as the dissolution of mathematics as a coherent academic discipline.

Education and research in mathematics occur outside mathematics departments because the social significance of ``mathematics'' remains close to its etymological meaning of all ``quantitative science'' as the latter is understood today. Of the 29,506 doctorates awarded annually by the 96 United States universities with very high research activity, 58\% were in social science, physical science, technology, engineering, or mathematics \citep {Carnegie2010}, all fields whose research is likely to involve mathematics. For example, the natural and social sciences formulate and test quantitative theories, which in effect are mathematical models, so progress often depends on mathematical prowess. Additionally, the social and health sciences gather voluminous survey and clinical data, which they analyze using mathematics developed in fields such as econometrics and epidemiology. Finally, the engineering and management sciences build models to mimic complex systems which they understand and may improve, through the models, using the mathematics of  game theory, operations research, and systems theory. Few subjects are so pervasive as mathematics in the intellectual life of the research university; consequently, expertise in mathematics is also ubiquitous.

\section {Research Universities in the United States}

The subject of this paper are United States research universities. By one assessment there are 96 institutions with very high research activity, and another 103 with high research activity \citep {McCormickZhao2005, Carnegie2010}. Nevertheless, while many institutions declare they support research and attempt to gain financial resources for it, only about 50 universities succeed in obtaining significant funding and perform the majority of notable research \citep [p.\ 58] {Calhoun2000}. Their effort is considerable. Surveys indicate that the faculty at such universities devote equal amounts of time on average to teaching and research, while the proportion of research time at other doctorate-granting institutions is very much less \citep [p.\ 330] {Finkelstein2001-in}. 

In the United States and in many other countries the smallest administrative units within universities are \textit {departments\/} \citep {Walvoord2000}. These units formulate curricula, teach the courses, originate the hiring of staff, and initiate the promotion of ladder staff.\footnote {\textit {Ladder\/} staff are faculty members who have been given  tenure or a promise to be considered.} The departments of large or research universities specialize in one academic subject. This tradition of specialization originated at some European universities in the 19th century and was partly responsible for promoting academic disciplines by creating job markets for faculty with specialized credentials.\footnote {\citet {Clark2006} amusingly describes the ``academic-charisma'' that universities seek in faculty members.} This symbiotic relationship between departments and disciplines is not entirely benign because the employment prospects of specialists implicitly depend on maintaining well demarcated disciplines.\footnote {A good example of the importance of the academic job market to disciplines is the \textit {Annual Survey of the Mathematical Sciences in the United States\/}, which in spite of the ambitious title is only an accounting of supply and demand for mathematics department faculty \citep {PhippsETAL2009}.} In this sense, disciplines are ``antithetical to a strong local intellectual community'' within the university \citep [74--75] {Calhoun2000}. These deep, sociological roots of specialization are perhaps overlooked when departments are criticized for ``inattention'' \citep [25] {Walvoord2000} to the interdisciplinary needs of undergraduate education \citep {Lattuca1994}. 

The emphasis on interdisciplinary undergraduate education is a curricular feature that especially distinguishes universities in the United States from universities in many other countries. The need for interdisciplinary eduction arises because entering students are assumed not to have completed their general educations \citep [p.\ 51]  {Calhoun2000}, and the university is responsible for providing that education. Thus, in the nominally 4 year curricula for bachelor's degrees, the \textit {lower division\/} years 1--2 are devoted to introductory courses from many fields.  Only during the \textit {upper division\/} years 3--4 do students concentrate on advanced subjects usually from one field, or what is often the same thing, from one department.

\section {Method}

\subsection {Empirical Methodology}

The methodology of this paper differs from previous sociological studies of mathematics in several respects. Most significant is the examination of mathematics throughout the university rather than only in eponymous departments. All previous studies implicitly equate mathematics with the intellectual content of papers and syllabi written by the faculty of mathematics departments. It is novel in the extreme to investigate where else mathematics expertise can be found. Further, only quantitative data is gathered to reveal the actual roles of participants rather than perceptions of those roles from opinion survey data. Mathematics activity is observed through databases for university courses, enrollments, and publications. Comments on the straightforward procedures that were employed to gather each type of data are included with the results. 

\subsection {Institutional Sample}
\label {sec:sample}

The data used throughout this paper are from 50 institutions that are profiled in \mytables {\ref {tab:overview} and \ref {tab:personnel}} of the appendix. These institutions have the most highly rated mathematics departments in the following sense. The National Research Council produced the most authoritative recent ranking of mathematics research-doctorate programs in United States universities \citep [appendix table H-4] {Goldberger1995}. Each department received a composite, numerical rating between 0.00 and 5.00. The American Mathematical Society groups departments with similar ratings to report annual statistics on employment in mathematics departments \citep {AMSsurvey, PhippsETAL2009}. Group I consists of the 48 most highly ranked departments, with ratings from 3.00 to 5.00. All but two of the 48 are from the 96 universities with very high research activity in the Carnegie classification, and the remaining 2 are from the 103 universities with high activity. The group I cadre thus provides a reasonably objective and quite broad sample of mathematics research departments. Two additional institutions are included because they provide some rarely available data that is used in section \ref {sec:coursework}.

\subsection {Identifying Mathematics Research}

This investigation requires a method to recognize, quantify, and classify mathematics research. Fortunately, that burden of subjectivity is assumed by two comprehensive surveys of research publications. Unlike the social sciences with which readers may be more familiar, many of the physical sciences index their research publications. These indexes include all peer-reviewed publications from any source that skilled editors judge to have contributed to the field. Mathematics has two indexing agencies sponsored respectively by the American and European Mathematical Societies. 

Both \citet {MathematicalReviews} and \citet {Zentralblatt} maintain databases of papers indexed by mathematics subject. Currently, \textit {Zentralblatt\/} surveys 3,500 journals and 1,100 serials for mathematics papers and indexes approximately 100,000 new peer-reviewed publications per year. 
Since 1971 \textit {Zentralblatt\/} has collected 2.39 million publications all of which are indexed by the Mathematical Subject Classification (MSC). This subject index consists of alphanumeric codes beginning with 2-digit numbers that reflect the coarsest level of differentiation, that is, the major branches of mathematics. Sixty-three of these 2-digit numbers are assigned.\footnote {The \citet {MSC2010} posts the MSC on the world wide web. \citet {Fairweather2009} discuss the latest revision. The MSC changes very infrequently at the 2-digit level, most recently in 2000 when class 04 was combined into 03.}

\begin {figure} 
\centering
\includegraphics [scale=1] {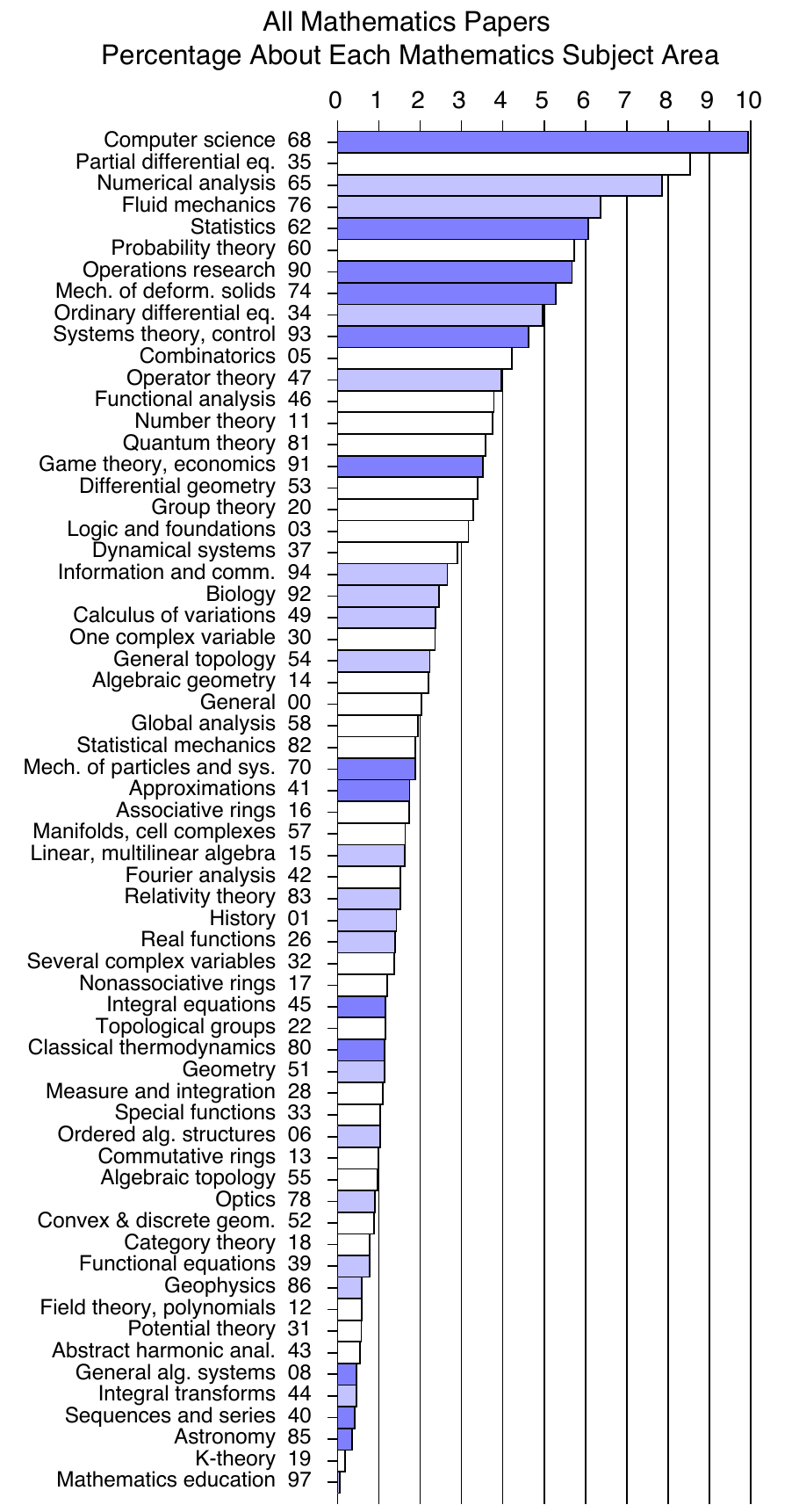}
\caption {A census of mathematics in the form of the percent of all mathematics papers addressing a particular Mathematics Subject Classification. The data were gathered from the Zentralblatt MATH online database in spring 2010. The percentages sum to 155.7 because many papers have more than one subject classification. Subjects are shaded as in Figure \ref {fig:emphasis}.}
\label {fig:census}
\end {figure}

For the purposes of this paper, mathematics research is quantified by the numbers of peer-reviewed papers indexed by \textit {Zentralblatt\/}. This data can be used to make a census of mathematics papers in the form of the percent of papers that address each major subject class (\myfigure {\ref {fig:census}}). For example, when the data were gathered in spring 2010, the \textit {Zentralblatt\/} database had 86,517 papers with an index code in the ``quantum theory'' classification, 81.\footnote {In comparison, the  \citet {Scitation2009} database had 216,747 papers on ``quantum'' subjects, thereby illustrating the selectivity of \textit {Zentralblatt\/} editors in choosing only mathematically relevant papers. \citet {Roth2005} describes publication databases for the physical sciences.} As noted, 2.39 million papers in the database had been assigned codes. Thus, approximately 86,517$ / (2.39 \times 10^6)$ or 3.62 percent of all mathematics papers contribute to the mathematics of quantum theory. This value is recorded in Figure \ref {fig:census}. Since a paper may have multiple subjects, the percentages in the figure sum to 155.7. The 10 subjects listed atop Figure \ref {fig:census} account for half the 2.39 million peer-reviewed papers. 

This census of mathematical research reflects current interests because most of the indexed papers have been written in the recent past. The quantity of mathematics publications has consistently increased by approximately 2,000 papers per year in recent decades, with the result that 68\% of the papers indexed by \textit {Zentralblatt\/} are from 1990--2009. If the data in figure \ref {fig:census} were restricted to just the past two decades, then the primary change among the heavily published subjects is to increase the percent of papers about the mathematics of computer science because computer science is a comparatively new field. 

\section {Results}

\subsection {Mathematics Department Coursework}
\label {sec:coursework}

Undergraduate instruction is a natural starting point for an investigation of mathematics in universities. The numeraire for measuring university instruction is the student credit-hour, SCH. The acronym is sometimes misinterpreted, and the concept is better understood, as a student-contact hour. An SCH represents one student who is being ``credited,'' for purposes of earning a degree, with one hour per week of supervised study per academic term. The concept arose in the United States as a quantitative measure of educational achievement to replace comprehensive examinations; the latter persist in Europe and in graduate studies. \citet {Shedd2003a, Shedd2003b} describes the history of the credit-hour and the vagaries of its calculation. From the standpoint of the university, credit-hours or contact hours measure how much effort university employees devote to education each week. 

All universities compile credit-hour data for their own use and for accreditation purposes. Data are usually given to the public summarized by campus, college, or school, but only occasionally by department (equivalently, by subject) and in a very few cases in even more detail. The lack of public access is responsible for considering relatively few universities in this section; however, the results for them are so overwhelming that there is no evidence opposed to the conclusions.

\begin {figure} 
\centering
\includegraphics [scale=1] {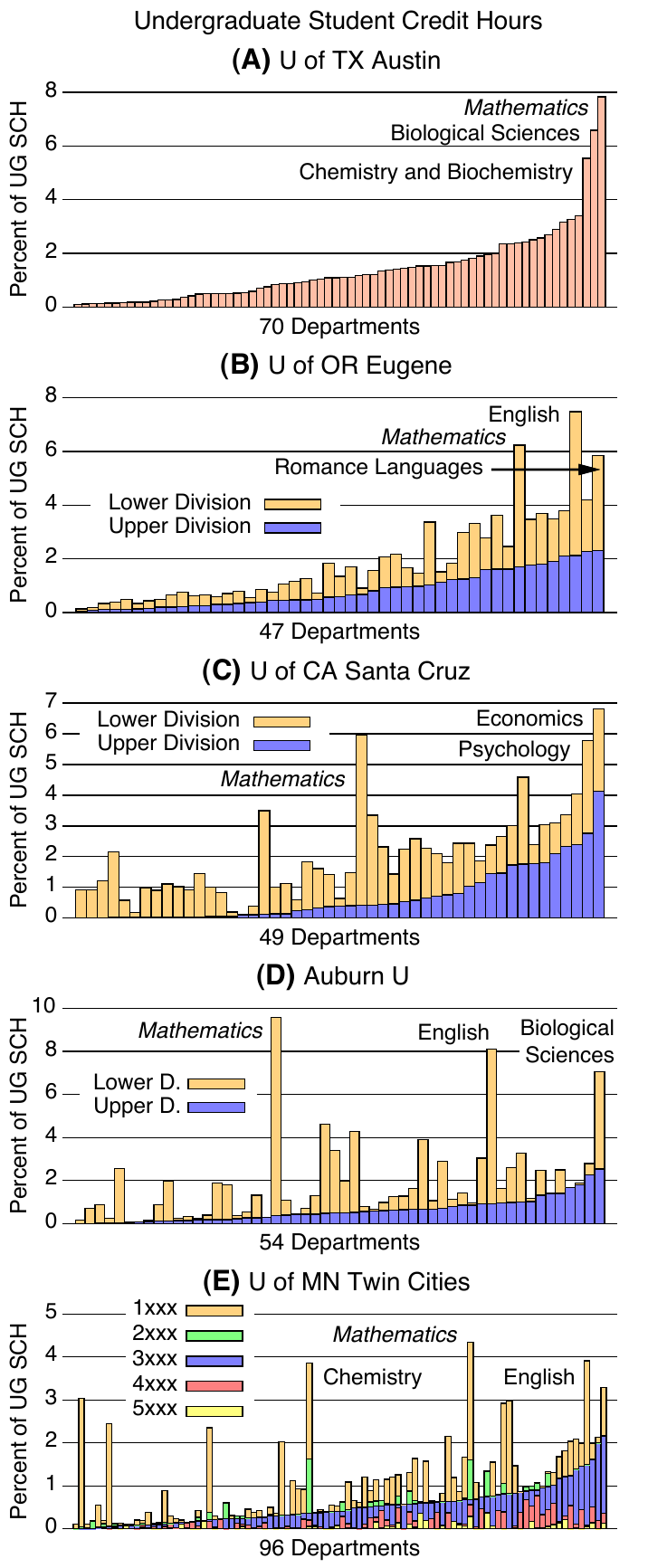}
\caption {Undergraduate student credit-hours by department for five universities. Departments are ordered by increasing upper division enrollment when the distinction is made. Only departments with values at least 0.1\% of credit-hours are shown. Data are for the fall 2008 semester except Minnesota data are for fall, spring, and summer 2005-2006. 
}
\vspace*{-2ex}
\label {fig:sch}
\end {figure}

The \citet [pp.\ 85--96] {UTXAustin2009} supplies credit-hour data by semester, department, and the degree level of the student. Undergraduate credit-hours for fall 2008 totaled 499,650 to which the Department of Mathematics contributed an astonishing 38,386 or 7.7\%. The mathematics department is the largest educational steward in Austin (\myfigure {\ref {fig:sch}A}). This rank seems to be held generally. Mathematics departments are among the two largest in terms of undergraduate credit-hours for each of the few universities providing detailed data to the public. Perhaps more remarkable is that, across all the universities considered, mathematics is the only department that consistently ranks among those with the most undergraduate credit hours.

More insight into mathematics coursework comes from universities that report data for more finely differentiated student levels. The \citet {UOREugene2009}, the \citet {UCASC2009}, and \citet {AuburnU2010} report credit-hours separately for lower division undergraduate students and for upper division students. These universities are comparably sized with respectively 236,302, 225,801, and 267,607 undergraduate student credit-hours in fall 2008. The mathematics departments supply 6\% of the undergraduate student credit-hours at the first two institutions and a very high 10\% at the third (\myfigures {\ref {fig:sch}B--D}). The more detailed data reveal most of these credit-hours are for lower division students. The upper to lower division ratio of SCH in mathematics courses is 0.37 : 1 in Eugene and a remarkably small 0.074 : 1 and 0.037 : 1 at Santa Cruz and Auburn, respectively.

Still greater insight comes from the \citet {UMN2006a} campus in the ``Twin Cities'' of Minneapolis and Saint Paul. This campus is comprehensive in the extreme, with professional schools of agriculture, business, engineering, law, and medicine. The university had 947,540 undergraduate student credit-hours in total for the fall, spring, and summer terms of the 2005--2006 academic year. Minnesota credit-hour data are reported by course numbers that reflect the year after matriculation in which students are anticipated to take the courses. Freshman courses are ``1xxx,'' sophomore ``2xxx,'' and so on. The \citet {UMN2006b} indicates this course numbering scheme is informally honored, so the numbers provide a reasonable indicator of course sophistication. Interpreting lower division courses as those numbered ``1xxx'' and ``2xxx,'' the upper to lower division ratio of credit-hours for the School of Mathematics is 0.19 : 1, consistent with the other schools. The observation enabled by this data is, mathematics departments provide by far the majority of their credit-hours in freshman courses (\myfigure {\ref {fig:sch}E}).

\subsection {Upper Division Mathematics Education}

The conceit, that science and technology are impossible without sophisticated mathematics, can be reconciled with Figures \ref {fig:sch}B--\ref {fig:sch}E (which show that mathematics departments mostly teach lower division courses) only if substantial mathematics education occurs outside mathematics departments. The extent of this teaching is difficult to quantify because course catalogs do not itemize the pertinent courses under the rubric ``mathematics.'' Consequently, class catalogs, course descriptions, syllabi, and textbooks all have to be examined to identify courses with mathematical content.\footnote {\citet [chap.\ 2] {Clark2006} refers to the academic catalog as ``the single most condensed academic document, the royal road to the academic subconscious.''} These documents at the \citet {UTclasses} are easily perused, and the coursework at this institution has already been examined in section \ref {sec:coursework}, so the very large student body of this research university is taken as representative.\footnote {Only three institutions with very high research activity have over 50,000 students \citep {Carnegie2010}.} All upper division courses offered in the fall and spring semesters of the 2008--2009 academic year are canvassed.

Mathematics courses are identified by matching them to specific items in the Mathematics Subject Classification. This approach assures a rigorous selection process. For example, the Department of Finance offered FIN 357, Business Finance. One professor \citep {Rao2008} assigned chapter 4 of the text \citep* {RossETAL2008} ``Discounted Cash Flow Valuation,'' and he scheduled lectures on net present value. This course therefore taught subjects under Mathematics Subject Classification 91B28, ``Finance, portfolios, investment,'' or most broadly under subject 91, ``Game theory, economics, social and behavioral sciences.'' There were 15 sections of FIN 357 in the fall and 8 in the spring, each meeting 3 hours per week, totaling 69 weekly class-hours devoted to subject area 91 (game theory, economics). 

For further example, 80\% of the 203 upper division class-hours in Department of Computer Sciences courses have mathematics content. Some courses are associated with more than one subject area. CS 320N, Practical Linear Algebra I \citep {vandeGeijn2009}, teaches how ``to attain high performance on current cache-based architectures,'' which is Mathematics Subject Classification 68N19, ``Other programming techniques.'' The course textbook \citep {Strang2006} discusses topics in Mathematics Subject Classification 15, ``Linear and multilinear algebra, matrix theory.'' Consequently, the 6 class-hours for this course are divided between subject classes 15 (linear algebra) and 68 (computer science). 

Altogether 252 upper division courses were identified with 667 sections meeting 2,197 hours per week during the fall or spring semesters. The distribution of these upper division class-hours with respect to the branches of mathematics and the departments offering the courses gives a reasonably accurate and very novel quantification of upper division instruction in mathematics subjects (\myfigure {\ref {fig:courses}}). Several areas of mathematics individually contribute 4\% or more of class-hours, as follows:

\begin {figure}
\centering
\includegraphics [scale=1] {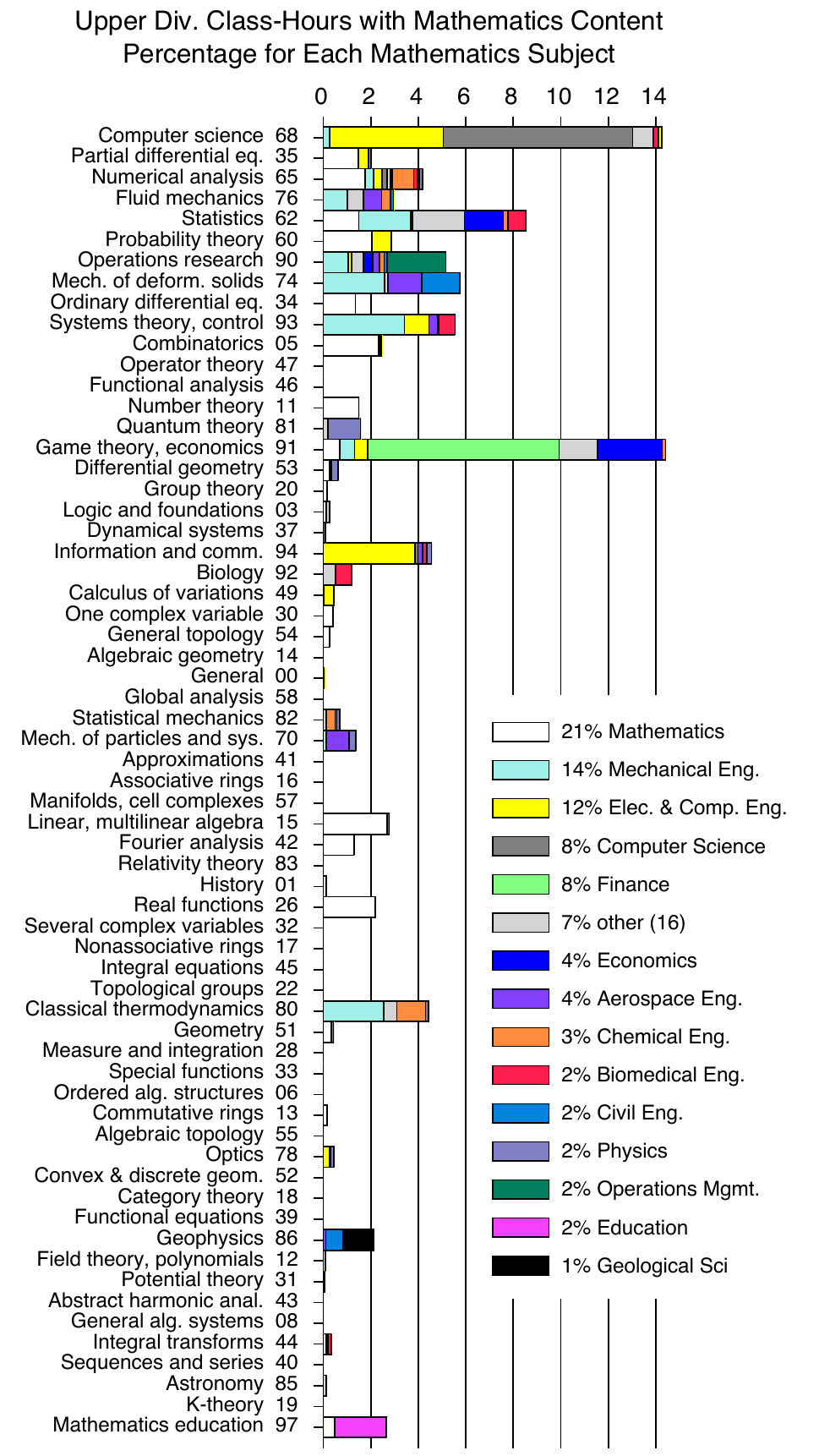}
\caption {A census of mathematics instruction in upper division courses at the University of Texas at Austin in fall 2008 and spring 2009. Student credit hours for courses with mathematics content are reported by the percentage that address each mathematics subject. The inset shows the distribution by academic department.}
\label {fig:courses}
\end {figure}

\begin {itemize}
\item \textbf {Computer science}, mathematics subject code 68, accounts for about 14\% of the upper division class-hours. The Departments of Computer Sciences and of Electrical and Computer Engineering provide most of the instruction in the mathematics of computer science; their courses overlap considerably as a result. 

\item \textbf {Economics}, code 91, also accounts for about 14\% of mathematics class-hours, but in several departments. The Department of Finance teaches how to value monetary instruments in business finance courses such as the one discussed above, FIN 357. Many engineering departments teach similar mathematics in courses about the financial management of engineering projects. The Department of Economics teaches the broadest range of subjects including game theory, ECO 354K, and econometrics, ECO 341K. The latter is one of those courses that address two areas of mathematics, 91 and 62. 

\item \textbf {Statistics}, 62, accounts for 8\% of the class-hours, the third largest. If the University of Texas at Austin had a statistics department, it likely would inherit the class-hours of the Department of Mathematics, which are only a quarter of those offered. The subject is also taught in many departments of engineering and social science. The campus does not have a medical school, else courses in biostatistics and epidemiology might well increase the class hours in statistics to rival computer science and economics. 

\item If \textbf {mathematical physics} were considered a separate subject class to consolidate the mathematics of fluid and solid mechanics 76 and 74, statistical and orbital mechanics 70, and thermodynamics 80, then it would account for over 14\% of class-hours. The Department of Mechanical Engineering teaches half the class-hours for these subjects.

\item \textbf {Control theory}, 93, is mostly taught in the Department of Mechanical Engineering, although at other universities it may also be found in aerospace engineering.

\item The Department of Electrical and Computer Engineering alone teaches the mathematics of \textbf {information and communication theory}, 94. 

\item The Department of Mathematics teaches half the class hours in \textbf {numerical analysis}, 65, with the rest taught by a variety of engineering departments. 

\item \textbf {Operations research}, 90, is taught in the Department of Information, Risk and Operations Management with the rest again taught by various engineering departments. 
\end {itemize}
This data has been gathered from a single university but the conclusion is broadly supported, that much mathematics education occurs outside mathematics departments. Since mathematics is used throughout the sciences, and since mathematics departments primarily teach lower division courses (as seen in section \ref {sec:coursework}), the necessary instruction in specialized or advanced mathematics must occur elsewhere.

The 2,197 class hours weekly in mathematics give a rough estimate of the student credit-hours. The university does not provide registration data for individual courses; however, Texas statutes require universities to ``offer only such courses and teach such classes as are economically justified'' \citep [title 3 ``education code'' \S51.403(d)] {TexasLegislature}. The \citet {TexasLA2009} interprets this law to mean undergraduate courses should have at least 10 students. Upper division courses are not anticipated to be large, so it is reasonable to suppose an average section size of 25. The 2,197 class-hours thus suggest 54,925 weekly student-contact hours in upper division mathematics instruction in total for both semesters. Recall, the university has approximately one-half million undergraduate student credit-hours per semester. 

With at most only about 5\% of student-contact hours devoted to upper division mathematics, the methods of instruction must be particularly effective. Even though most courses may be taught in conventional lectures, syllabi that contain a mixture of mathematics and other subjects are unconventional in their treatment of mathematics. They introduce the mathematics early, so the remainder of the course effectively becomes a practicum for the mathematics. For example, PHY 333, Modern Optics \citep {Keto2009}, successively treats geometric optics, vectors and matrices in the Jones calculus, and Fourier transforms, each followed by studies of optics using this mathematics. Mathematics departments have recognized the value of this type of instruction where it is generally called ``experience in model building.'' Among the 2,197 section hours, the mathematics department has 3 for M 474M, Mathematical Modeling in Science and Engineering. Most of the courses outside the mathematics department teach mathematics in this way.

\subsection {Research Emphasis of Mathematics Departments}

Having obtained data on instruction, the discussion turns to research activity. The 48 leading mathematics departments have 2,108 ladder staff with primary appointments in them as of spring 2010. Their research papers were obtained from the \textit {Zentralblatt\/} database by subject classification and by university. This method avoids double-counting jointly-written papers when all authors belong to the same institution: 63,874 papers were obtained. These peer-reviewed papers are all research of all ladder faculty at all departments that the American Mathematical Society regards as leading research in mathematics; therefore, there is no sampling error in the following discussion.

\begin {figure} 
\centering
\includegraphics [scale=1] {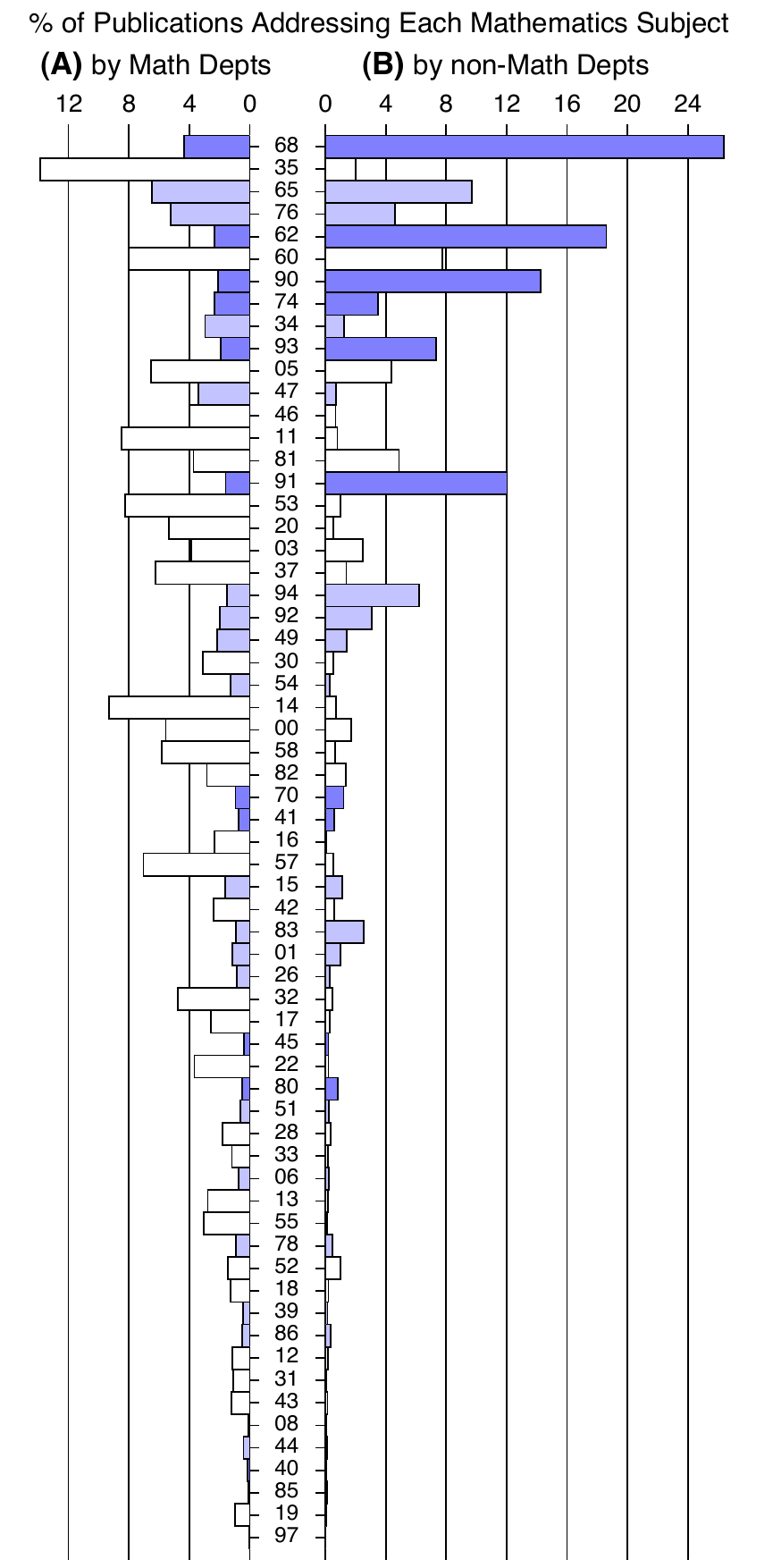}
\caption {Of the papers written (A) by mathematics departments, and (B) by non-mathematics departments, the percent that address each mathematics subject. Percentages sum to 184.0 and 155.7, respectively, because many papers have more than one subject. Subjects are shaded as in Figure \ref {fig:emphasis}.}
\label {fig:ext-and-int}
\end {figure}

The 63,874 papers have the subject distribution shown in \myfigure {\ref {fig:ext-and-int}A}. Comparing Figures \ref {fig:census} and \ref {fig:ext-and-int}A indicates the research interests of leading mathematics departments are unrepresentative of all mathematics. For example, in Figure \ref {fig:census}, about 1.6\% of all mathematics papers are about manifolds and cell complexes 57, while in Figure \ref {fig:ext-and-int}A, mathematics departments write 7.0\% of their papers about this subject. Perhaps more significant, 9.9\% of all mathematics papers in Figure \ref {fig:census} are about computer science 68, the most of any subject, while in Figure \ref {fig:ext-and-int}A only 4.3\% of the papers written in mathematics departments are about the mathematics of computer science. 

\begin {figure}
\centering
\includegraphics [scale=1] {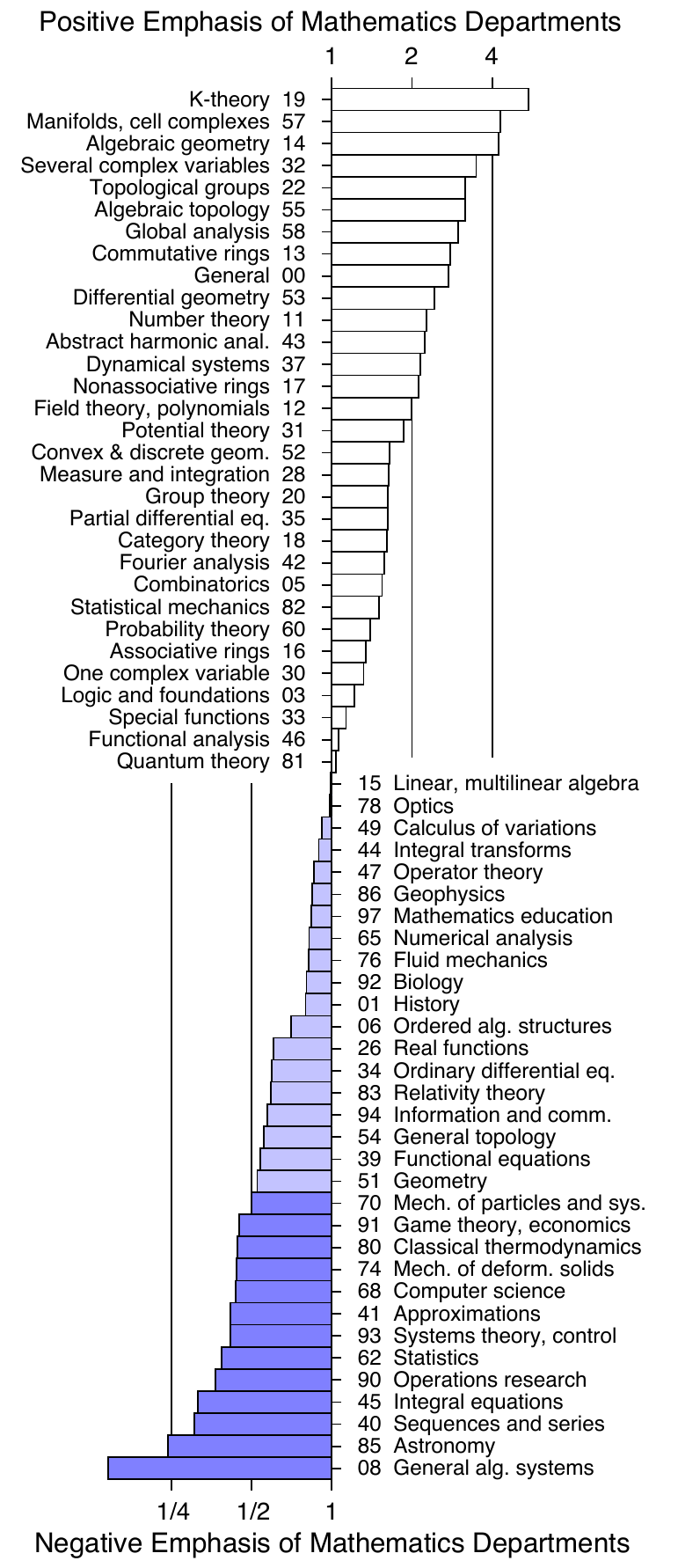}
\caption {Publication emphasis of mathematics departments, the ratio of the fraction of department publications that are about a given subject to the fraction of publications in all of mathematics about the subject. Note the logarithmic scale. Shading is repeated in Figures \ref {fig:census} and \ref {fig:ext-and-int}.}
\label {fig:emphasis}
\vspace*{-3ex}
\end {figure}

The research emphasis of a group of mathematicians may be defined as the ratio of the fraction of their papers on a subject to the fraction of all mathematics papers on the subject. This quantity provides a criterion to identify subjects in which mathematics departments have more or less expertise. It is convenient to represent the ratio as a base $2$ logarithm, so an emphasis in favor has a positive value and an emphasis against has a negative value. \myfigure {\ref {fig:emphasis}} reports the research emphasis of mathematics department ladder staff. Mathematics departments have a strong negative emphasis (less than $2^{-1}$) against 13 subjects that are darkly shaded in Figures \ref {fig:census}, \ref {fig:ext-and-int}, and \ref {fig:emphasis}. They have a moderate negative emphasis (between $2^{-1}$ and $2^{\kern0.05em 0}$) against 19 subjects that are lightly shaded in the figures. Remarkably, these figures show mathematics departments have a negative research emphasis for 8 of the 10 areas of mathematics at the top of Figure \ref {fig:census} about which the most research is conducted.

\subsection {Mathematics Research Outside the Mathematics Department}

Mathematics departments de-emphasize research in the most heavily published mathematics subjects (comparing Figure \ref {fig:emphasis} with Figure \ref {fig:census}). Consequently others must write many of the papers. In order to identify the authors, this section examines departments other than mathematics at the same 48 universities with the leading mathematics departments. \mytable {\ref {tab:departments}} in the appendix lists the selected departments. They have been chosen to be among those that Figure \ref {fig:courses} shows provide much of the upper division mathematics instruction at the University of Texas. These non-mathematics departments had 13,701 ladder staff with primary appointments in them in spring 2010.\footnote {Only faculty members with primary appointments outside mathematics departments are considered in case of joint appointments. When a university has two mathematics departments then only faculty outside both are considered.} Like the investigation of mathematics departments, papers in the \textit {Zentralblatt\/} database are gathered by subject classification for the whole sample at each university to avoid double-counting jointly-written papers.  The subject distribution of the 75,125 papers that were written by members of the chosen non-mathematics departments indicates that the mathematics research interests of authors outside mathematics departments are largely complementary to those inside (\myfigure {\ref {fig:ext-and-int}B}). 

\begin {figure} 
\centering
\includegraphics [scale=1] {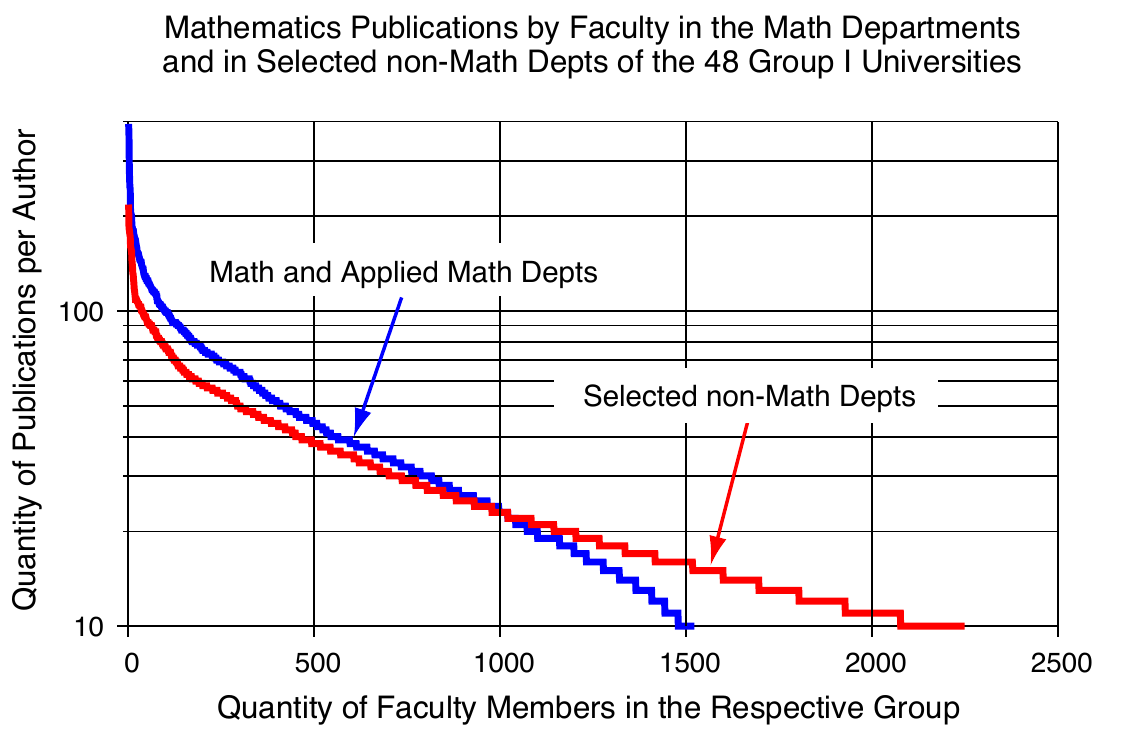}
\caption {Quantity of peer-reviewed mathematics publications by each ladder faculty member of mathematics, or of selected non-mathematics departments, at the 48 Group I universities with the most highly rated mathematics departments. Only authors with 10 or more mathematics publications are shown. Table \ref {tab:departments} in the appendix lists the selected non-mathematics departments.}
\label {fig:productivityBYauthor}
\end {figure}

A reasonable measure of expertise in any subject is having written 10 peer-reviewed publications. The ratio of such authors outside and inside the mathematics departments of the 48 universities is at least $2249 : 1522 = 1.47$ (\myfigure {\ref {fig:productivityBYauthor}}). Thus research universities have more mathematically expert faculty outside their mathematics departments than inside. The ratio can only increase by including more departments of natural or social science in the selection of non-mathematics departments. With the present sample the mathematics departments have more authors of mathematics papers only for authors of 22 or more papers. Writing so many mathematics papers is a remarkable achievement for faculty outside mathematics departments because they also may write papers in their primary fields of study. 

\begin {figure} 
\centering
\includegraphics [scale=1] {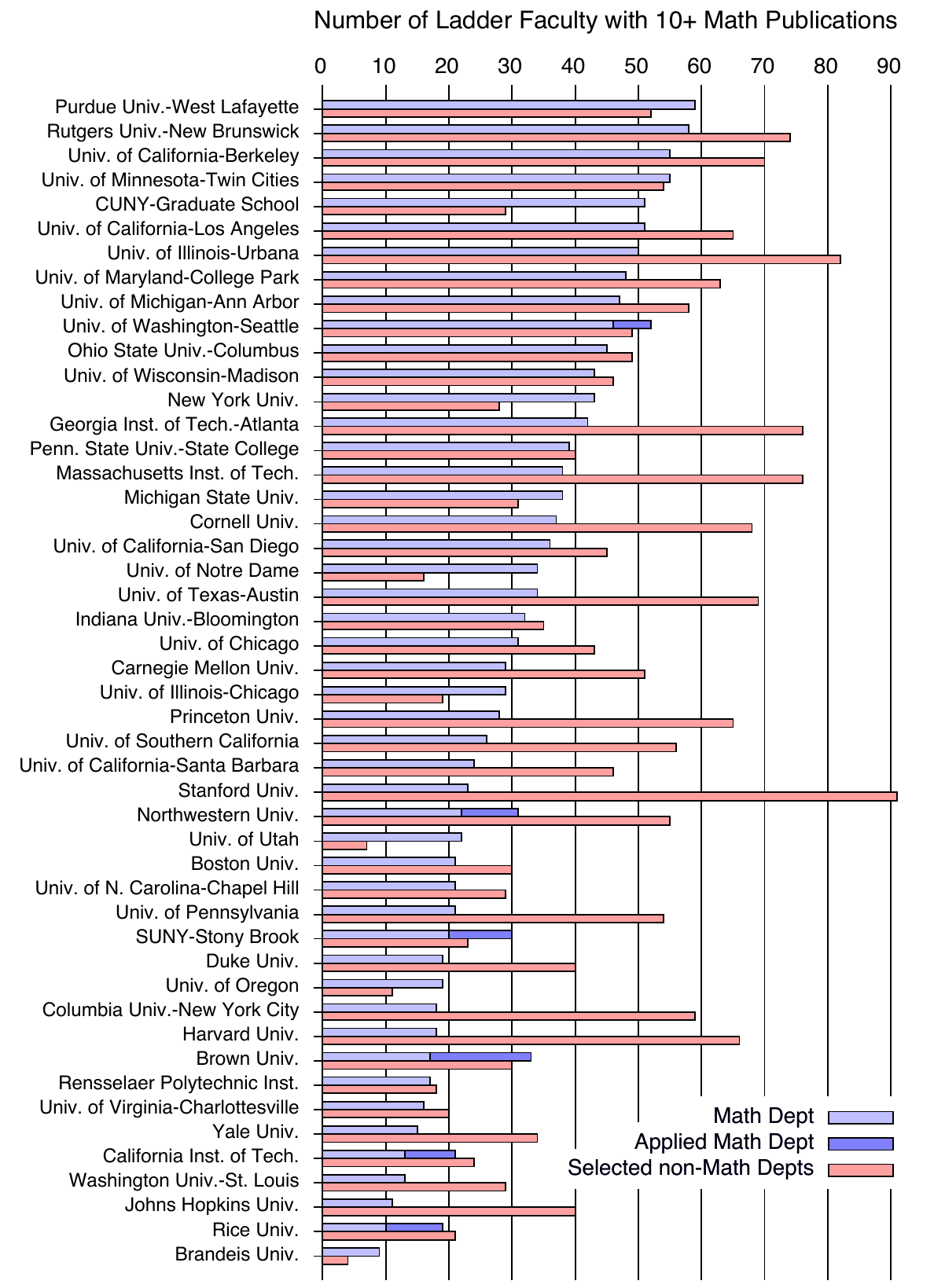}
\caption {Quantity of ladder faculty with 10 or more peer-reviewed mathematics publications in the mathematics departments, or in selected non-mathematics departments, of the 48 Group I universities with the most highly rated mathematics departments. Universities are listed in order of number of faculty in their mathematics department. A few universities have a second (applied) mathematics department whose faculty are shown separately. Table \ref {tab:departments} in the appendix lists the selected non-mathematics departments.}
\label {fig:productivityBYdept}
\end {figure}

Faculty who conduct research in mathematics are equally productive inside and outside mathematics departments. Based on the criterion of having written 10 or more mathematics papers, the contingent outside mathematics department is larger at 38 of the 48 institutions (\myfigure {\ref {fig:productivityBYdept}}). This group is smaller at some institutions such as Brandeis and New York University that lack extensive collections of professional schools. A characteristic of excellence in comprehensive research universities may be that mathematics research skills are distributed throughout the faculty.

\section {Discussion}

\subsection {Intensive and Extensive Nomenclature}
\label {sec:nomenclature}

The aforementioned mathematics faculties inside and outside mathematics departments are not discussed in the sociology of science or of higher education, so some new naming convention is needed to refer to them. Borrowing terminology from the physical sciences, the ``intensive'' mathematics faculty are members of mathematics departments. The ``extensive'' faculty are tenured and tenure-track staff who teach mathematics or write mathematics papers outside mathematics departments. The analogy with physical properties is, the quantity of extensive faculty varies in proportion to the size of the whole natural and social science faculty. 

The extensive faculty are not necessarily mathematicians by training or preference. Nevertheless, they perform as such, by teaching upper division courses whose syllabi include mathematics subjects not taught elsewhere, and by writing the majority of papers that bibliographic services regard as having mathematics content. The observation, that significant research activity for any field occurs apart from the academic department nominally associated with the field, appears to be unprecedented. 

The extensive faculty do not coincide just with those faculty outside mathematics who conduct research on mathematics subjects. Many faculty may teach mathematics subjects that they know but on which they do not publish because teaching duties typically are broader than research interests. Nevertheless, the results of this study show the extensive faculty contains sizeable numbers of individuals who do publish significantly in mathematics (Figures \ref {fig:productivityBYauthor}, \ref {fig:productivityBYdept}), and moreover these faculty produce the bulk of research about certain mathematics subjects (Figure \ref {fig:ext-and-int}B).

\subsection {Convergence or Divergence of Teaching and Research}

The results of this study reveal a significant dichotomy between the extensive and intensive mathematics faculties. The extensive faculty write papers on the mathematics s ubjects they teach, whereas mathematics departments write papers on subjects that are largely separate from their undergraduate teaching. 

A meaningful data reduction begins by selecting the mathematics subject areas that are most important in the teaching or publications of the extensive or intensive faculties. Publication subjects are chosen by applying a 4\% threshold to Figures \ref {fig:ext-and-int}A and \ref {fig:ext-and-int}B; that is, a subject must be addressed in either 4\% of papers written by the extensive faculty, or separately 4\% of papers written by the intensive faculty. Teaching subjects are chosen similarly from a 3\% threshold on the class-hours in Figure \ref {fig:courses}. This figure reports consolidated percentages that must be adjusted for purposes of selection to reflect class-hours of the extensive or intensive faculties alone. This selection of subjects paints an objective picture of what mathematics is taught or published in quantity and by whom (\mytable {\ref {tab:most}}). The thresholds for inclusion are quite low so any subject omitted is genuinely secondary compared to others. 

\begin {table} [h!]
\renewcommand {\baselinestretch} {1.0}
\caption {Mathematics subjects as interest in them is apportioned among the extensive or intensive faculties. Teaching subjects are determined by applying a 3\% threshold in Figure \ref {fig:courses} to the courses taught by the respective group. Research subjects are those that attain 4\% in Figures \ref {fig:ext-and-int} or \ref {fig:census}. Overall interest is determined by applying the 3\% and 4\% thresholds to Figures \ref {fig:courses} and \ref {fig:census}, respectively.}
\label {tab:most}
\begin {center}
\medskip
\newcommand {\bld} [1] {\textbf {#1}}
\newcommand {\ultxt} [1] {\underline {\smash {#1}}}
\newcommand {\ulbld} [1] {\underline {\smash {\textbf {#1}}}}
\renewcommand {\arraystretch} {1}
\scriptsize
\setlength {\tabcolsep} {0.33em}
\begin {tabular} {r c | c c | c c | c c |}
&&
\multicolumn {2} {c |} {Math.}&
\multicolumn {2} {c |} {non-Math.}&
\multicolumn {2} {c |} {Overall}\\
\cline {3-8}
\multicolumn {1} {c} {\first Subject}& \multicolumn {1} {r |} {\hspace*{-1em}Code}&
\multicolumn {1} {c |} {tea.}& \multicolumn {1} {c |} {res.}&
\multicolumn {1} {c |} {tea.}& \multicolumn {1} {c |} {res.}&
\multicolumn {1} {c |} {tea.}& \multicolumn {1} {c |} {res.}\\
\hline \first 
Combinatorics&05&\chk&\chk&&\chk&&\chk\\
Number theory&11&\chk&\chk&&&&\\
Algebraic geometry&14&&\chk&&&&\\
\hline \first 
Linear, multilinear algebra&15&\chk&&&&&\\
Group theory&20&&\chk&&&&\\
Real functions&26&\chk&&&&&\\
\hline \first 
Several complex variables&32&&\chk&&&&\\
Ordinary differential eq.&34&\chk&&&&&\chk\\
Partial differential eq.&35&\chk&\chk&&&&\chk\\
\hline \first 
Dynamical systems&37&&\chk&&&&\\
Fourier analysis&42&\chk&&&&&\\
Functional analysis&46&&\chk&&&&\\
\hline \first 
Differential geometry&53&&\chk&&&&\\
Manifolds, cell complexes&57&&\chk&&&&\\
Global analysis&58&&\chk&&&&\\
\hline \first 
Probability theory&60&\chk&\chk&&\chk&&\chk\\
Statistics&62&\chk&&\chk&\chk&\chk&\chk\\
Numerical analysis&65&\chk&\chk&\chk&\chk&\chk&\chk\\
\hline \first 
Computer science&68&&\chk&\chk&\chk&\chk&\chk\\
Mech. of deform. solids&74&&&\chk&&\chk&\chk\\
Fluid mechanics&76&&\chk&\chk&\chk&&\chk\\
\hline \first 
Classical thermodynamics&80&&&\chk&&\chk&\\
Quantum theory&81&&&&\chk&&\\
Operations research&90&&&\chk&\chk&\chk&\chk\\
\hline \first 
Game theory, economics&91&\chk&&\chk&\chk&\chk&\\
Systems theory, control&93&&&\chk&\chk&\chk&\chk\\
Information and comm.&94&&&\chk&\chk&\chk&\\
\hline 
\end {tabular}
\end {center}
\vspace*{0.1in}
\end {table}

\begin{figure} [b]
\begin{center}
\includegraphics [scale=1] {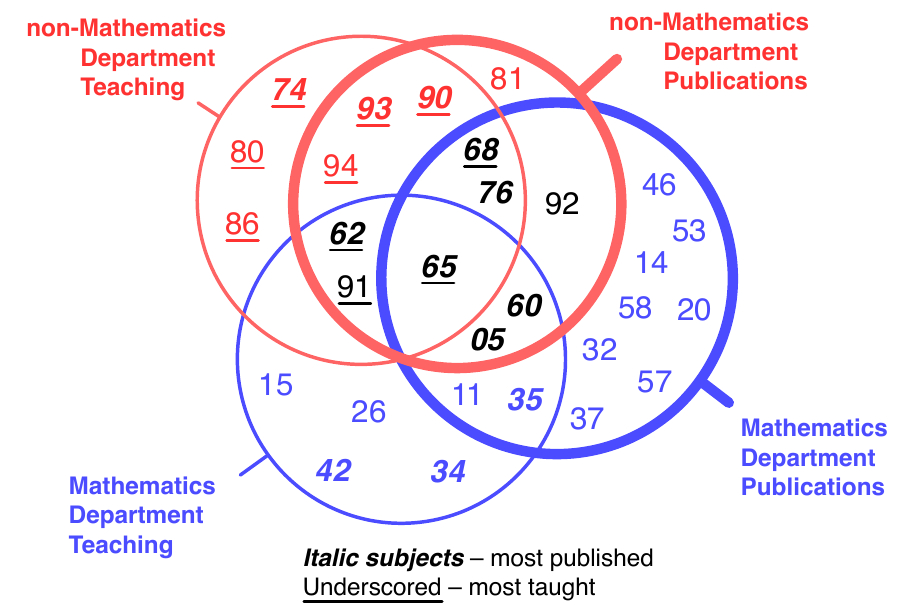}
\end {center}
\caption {Venn diagram of mathematics subjects that are prominent in publication or undergraduate teaching of either the extensive or the intensive mathematics faculties. Table \ref {tab:most} identifies the subjects. The text explains the selection criteria.}
\label {fig:venn}
\end {figure}

The subjects of interest to the intensive and extensive mathematics faculties can be arranged in a Venn diagram (\myfigure {\ref {fig:venn}}). Mathematics departments have only 5 subjects common to their publications and undergraduate teaching, and another 11 subjects on which they publish without appreciable instruction. In marked contrast, the non-mathematics departments have 8 mathematics subjects common to both their publications and undergraduate courses. This extensive faculty publishes on 8 of the 11 most heavily published mathematics subjects, it teaches all 9 most heavily taught subjects, and it is the primary teacher for 7 of them. 

\subsection {Mathematics as Meta-discipline}
\label {sec:meta}

Meta-disciplines have been broadly characterized as fields that encompass the representation of knowledge in other fields. They relate to other disciplines in a purely formal way: meta-disciplines have or contribute process-knowledge as opposed to domain-knowledge.\footnote {Information science has been called a meta-discipline in this sense \citep {Bates1999}. Fields that draw expertise from more than one other field have also been called meta-disciplines, for example see \citet {Mihelcic2010}. In the latter sense some fields may be meta-disciplines only temporarily while they are nascent.} This characterization would appear to be appropriate for mathematics.

If mathematics occurs in another field, it should not be presumed that only the application of a mathematical meta-process has transpired. To the contrary, often the relevant mathematics originated in the other field. For example, C.\ F.\ \citet {Gauss1809-HM} discovered the method of least squares, or regression analysis, during his work in astronomy and geodesy.\footnote {\citet {Stigler1986} describes this history, and \citet {Hald2007} derives the theory in the historical context. Gauss, whose career coincided with the formation of academic departments, is remembered for his mathematical and physical discoveries. He taught mostly astronomy and geodesy \citep [pp.\ 405--410] {Dunnington2004-HM}.} Today regression analysis is taught in mathematics or statistics departments (and others), yet Gauss's discovery was made while he was immersed in the physical sciences. Since discovery is the \textit {raison d'{\^e}tre\/} for research universities, it is reckless to ignore how mathematics discoveries have often been made. 

United States research universities have no consistent organizational structure for meta-disciplines such as mathematics whose expertise spreads across many fields. Considering just academic departments that involve considerable mathematics knowledge, the 48 group I universities (section \ref {sec:sample}) have disparate departments (Table \ref {tab:departments}) often in different colleges of the university. For example, operations research may be found in one or both of the business and engineering colleges. The mathematics of computer science is the purview of one or two departments variously called ``computer science,'' ``computer and information science'' and ``electrical and computer engineering'' that are found in engineering  or in science colleges. Separating computer science from mathematics has been described as detrimental \citep {Kowalik2006}. One may conjecture that the same applies to other subjects as well.

Statistics is a specific organizational challenge that all universities faced simultaneously and for which there is some historical narrative. Universities were urged to create statistics departments in the middle half of the 20th century because mathematics departments failed to add statisticians to their faculties in proportion to the growing demand for statisticians in government and industry. \citet {Hotelling1940} alleged that mathematics departments viewed teaching statistics as service work that did not require a commitment to research, thereby allowing new positions to be filled by specialists in other mathematics, thus decreasing the acuity of the statistics instruction.\footnote {If Hotelling is accurate, then the treatment of statistics in mathematics departments is a good example of intra-disciplinary considerations trumping those of the academic community \citep [pp.\ 74--75] {Calhoun2000}.} Whatever the case, administrations did not respond uniformly to their ``statistics problem,'' suggesting that they failed to resolve it. Some universities continue to teach statistics in their mathematics departments, while many have statistics departments now (Table \ref {tab:departments}), and all universities allow the subject to be taught by other departments as well \citep {Raftery2000}. These solutions may be reasonable except in so far as addressing what it means to maintain mathematics as a meta-discipline. Exhortations to be ``more interdisciplinary'' and ``more mathematical'' are a common stopgap.\footnote {See recommendations of the \citet [II, 2-20] {NSB2004} and several studies quoted therein.}

For several years, mathematics departments in the United States have had significant, precipitous declines in both undergraduate degrees and graduate enrollments, while at the same time, there has been essentially no decline in tenured faculty \citep {Bressoud2009, Kirwan2001}.\footnote {D.\ M.\ Bressoud is president of the Mathematical Association of America. W.\ E.\ Kirwan is chancellor of the university system of Maryland and a former mathematics department head.} As national economic conditions strain financial support for higher education, many institutions may look for economies in their mathematics programs. Faced with similar hardships a decade earlier, for example, one university proposed to employ only non-ladder staff for lower division mathematics instruction in the numerous service courses (to be supervised by natural and social science faculty), and to curtail the graduate mathematics program \citep {Jackson1996, Jackson1997}. Such proposals may engender vehement objections on the grounds that mathematics is a core academic discipline. However compelling that view may be, this paper shows the objections are not apropos the policy question, because the core of advanced undergraduate mathematics education occurs outside mathematics departments. The proper justification for maintaining mathematics does not appear to be that it is a core discipline, rather that it is a meta-discipline, but then in what form should it be maintained? 

\subsection {Philosophy and Sociology of Mathematics}

The failure to recognize the ramifications of mathematics as a meta-discipline, with a large ``extensive'' faculty outside mathematics departments, extends beyond university organizational structures. 

Studies of the sociology and philosophy of mathematics have focussed on the foundational issue concerning the epistemology of pure abstractions and its significance for education.\footnote {For example see the recent collection of articles edited by \citet {vanKerkhovevanBendegem2007}.} This approach amounts to restricting the view of mathematics to the intensive faculty and ignoring the extensive faculty. The pure abstractions in the research of the intensive faculty are doubly removed from reality because, in view of Figure \ref {fig:venn}, mathematics departments conduct research predominantly on subjects they do not teach to the educated public, who consist of university graduates with four-year degrees. Thus, the ``foundational question'' is a Gedankenversuch for philosophers and sociologists that is moot for society.\footnote {A \textit {Gedankenversuch\/} is a thought experiment in physics that is interesting from the standpoint of theory but is impossible to encounter in reality.} In contrast, the extensive mathematics faculty conduct research on tangibles whose meaning is undoubtable. For example, the professor of FIN 357 wrote a paper \citep {LeeRao1988} about asset valuation formulas, and the professor of CS 320N wrote a paper \citep {QuintanaGeijn2003} about devising algorithms for computers.\footnote {This particular professor of finance is not heavily published on mathematics subjects while the professor of computer science is.} A philosophical question about mathematics that would appear to be more relevant to society and education is the one asked by Eugene \citet {Wigner1960}: why is mathematics so successful in reasoning about reality? 

\subsection {K-12 Mathematics Education}

Twentieth century K-12 educational policy in the United States has been characterized as a continual disagreement between educators, mathematicians, and parents.\footnote {As described by \citet {Klein2002} and \citet {Schoenfeld2004}. In the United States, K-12 refers to kindergarten followed by 12 years of pre-university instruction.} How to prepare students to enter universities is a paramount concern, yet in view of the many departments that provide upper division mathematics instruction (Figure \ref {fig:courses}), participants in the K-12 debate may lack a comprehensive understanding of university mathematics. For example, K-12 educators acquire their own conception of mathematics from mathematics departments rather than from the majority of mathematics researchers who make up the extensive faculty. 
\begin {quotation}
\noindent
In an era of accountability and amidst brewing controversy concerning mathematics education's goals and emphasis, mathematicians and mathematics educators continued to work together to formulate recommendations for how teachers could develop the mathematical and pedagogical knowledge and skills they would need to help all students learn. \hfill  \mbox {--- \citet [p.\ 1268] {Ferrini-Mumdy2003}}
\end {quotation}

Natural and social science faculty recently insisted that statistics be included in curricular standards for university-bound students in spite of mathematics faculty who felt statistics was ``not necessarily a prerequisite for success'' \citep [p.\ 37] {Conley2003}. The suggestion that departments other than mathematics might set standards for K-12 education is considered to be novel \citep {Confrey2007}. Since the mathematics educations of many university students are completed by the extensive mathematics faculty, their view of the mathematics skills needed by university students would appear to be relevant to K-12 education.  

\subsection {National Research Education Policy}

At the turn of the last century, the United States \citet {NSF2000} began a grant program, VIGRE, for ``departments in the mathematical sciences to carry out high quality educational programs, at all levels, that are vertically integrated with the research activities of these departments.'' In short, the purpose was to improve university education in mathematics. Grants were awarded entirely to mathematics departments \citep {Durrett2002, Mackenzie2002}. These funds may have been misdirected because, as Figure \ref {fig:courses} indicates, 79\% of upper division mathematics student credit hours are taught outside mathematics departments to students who earn degrees in other departments. There is no systematic effort to address undergraduate university mathematics in a comprehensive manner. Referring to the extensive faculty, \citet {Steen2006} finds an exciting ``stealth curriculum that thrives outside the confining boundaries of college and university mathematics departments.''\footnote {L.\ A.\ Steen is a former president of the Mathematical Association of America.}

\section {Conclusion}

In summary, this paper shows that United States research universities have an ``intensive'' mathematics faculty inside a department of the same name, and an ``extensive'' mathematics faculty spread across other departments. The two faculties perform different roles in education and research. The intensive and extensive faculties teach lower or upper division students, respectively, and they conduct research in mathematics subjects either unrelated or relevant to upper division students, also respectively. United States research universities have not been successful in aligning the teaching and research responsibilities of their mathematics department faculties. In contrast, the mathematics research conducted by faculty outside the mathematics department is aligned with the vocational interests of undergraduate students. The existence of the extensive faculty has not been appreciated previously because an examination of the education literature indicates that the contributions of the faculty outside mathematics departments are not recognized in formulating national policies for mathematics education. 

This study suggests several questions for further investigation. One is whether other fields besides mathematics, such as other meta-disciplines, have similar ``two-tier'' faculties consisting of a formal, intensive faculty in one department and an informal, extensive faculty in other departments. Further is the question of how these separate faculties developed in the department-based organization of research universities and what they mean for the discipline-based organization of science. From the standpoint of educational administration it is relevant to study the  attitudes of the intensive and extensive faculty toward the other group and their perception of their own roles in the two-tier system. These attitudes evidently affect how university-wide curriculum is established. Coordination would appear to be through informal arrangements because direct administration would have to cut across the boundaries of departments and colleges which, ipso facto, are independent in formulating curricula. The utimate question is whether the two-tier approach is effective in promoting excellence in either research or education. As \citet [p.\ 48] {Calhoun2000} has written, the United States approach to education involves an interdependent system, but not necessarily one that reflects a rational design or that functions perfectly.

\section* {Acknowledgements}

I wish to thank the editor, Prof.\ C.\ Musselin, and the reviewers for comments and suggestions that greatly improved this paper.

\newpage

\appendix

\section* {Appendix: Institutions}

\begin {table} [h!]
\caption {Overview of  institutions discussed in this paper. ``R\&D'' is research and development. ``RU/(V)H'' is research university with (very) high research activity. ``\$M'' is one million dollars. ``U'grad'' is undergraduate. Data are from the \citet {Carnegie2010}.}
\label {tab:overview}
\begin {center}
\medskip
\scriptsize
\begin {tabular}
{l l r l l r r}
\multicolumn {3} {r |} {Science and Engineering R\&D Expenditures}&
\multicolumn {1} {r |} {\relax}&
\multicolumn {1} {c |} {Region of the}&
\multicolumn {2} {c |} {Enrollment}\\
\cline {6-7}
\first \textsc {institution}& 
\multicolumn {1} {r} {\hspace*{-8em}Carnegie Classification}&
\multicolumn {1} {r |} {\footnotesize \$M}&
\multicolumn {1} {c |} {Control}&
\multicolumn {1} {c |} {United States}&
\multicolumn {1} {c} {Total}&
\multicolumn {1} {r |} {\hspace*{-1em}U'grad.}\\
\hline \first
Auburn Univ.& RU/H& 117& Public& Southeast& 22,928& 78\%\\
Boston Univ.& RU/VH& 225& Private& New England& 29,596& 57\%\\
Brandeis Univ.& RU/VH& 51& Private& New England& 5,072& 63\%\\
Brown Univ.& RU/VH& 125& Private& New England& 8,004& 73\%\\
California Inst. of Tech.& RU/VH& 241& Private& Far West& 2,171& 41\%\\
Carnegie Mellon Univ.& RU/VH& 186& Private& Mid East& 9,803& 55\%\\
Columbia Univ.-New York City& RU/VH& 438& Private& Mid East& 21,648& 30\%\\
Cornell Univ.& RU/VH& 555& Private& Mid East& 19,518& 70\%\\
CUNY-Graduate School& RU/H& 3& Public& Mid East& 4,234& 0\%\\
Duke Univ.& RU/VH& 520& Private& Southeast& 12,770& 49\%\\
\hline \first
Georgia Inst. of Tech.-Atlanta& RU/VH& 364& Public& Southeast& 16,841& 65\%\\
Harvard Univ.& RU/VH& 409& Private& New England& 24,648& 32\%\\
Indiana Univ.-Bloomington& RU/VH& 135& Public& Great Lakes& 37,821& 75\%\\
Johns Hopkins Univ.& RU/VH& 637& Private& Mid East& 18,626& 28\%\\
Massachusetts Inst. of Tech.& RU/VH& 486& Private& New England& 10,320& 40\%\\
Michigan State Univ.& RU/VH& 321& Public& Great Lakes& 44,836& 73\%\\
New York Univ.& RU/VH& 242& Private& Mid East& 39,408& 48\%\\
Northwestern Univ.& RU/VH& 320& Private& Great Lakes& 17,747& 47\%\\
Ohio State Univ.-Columbus& RU/VH& 496& Public& Great Lakes& 50,995& 68\%\\
Penn.\ State Univ.-State College& RU/VH& 480& Public& Mid East& 41,289& 82\%\\
\hline \first
Princeton Univ.& RU/VH& 180& Private& Mid East& 6,708& 70\%\\
Purdue Univ.-West Lafayette& RU/VH& 309& Public& Great Lakes& 40,108& 76\%\\
Rensselaer Polytechnic Inst.& RU/VH& 51& Private& Mid East& 6,696& 73\%\\
Rice Univ.& RU/VH& 52& Private& Southwest& 4,855& 60\%\\
Rutgers Univ.-New Brunswick& RU/VH& 248& Public& Mid East& 34,696& 73\%\\
Stanford Univ.& RU/VH& 603& Private& Far West& 18,836& 35\%\\
SUNY-Stony Brook& RU/VH& 200& Public& Mid East& 21,685& 60\%\\
Univ. of California-Berkeley& RU/VH& 507& Public& Far West& 32,803& 67\%\\
Univ. of California-Los Angeles& RU/VH& 849& Public& Far West& 35,966& 68\%\\
Univ. of California-San Diego& RU/VH& 647& Public& Far West& 24,663& 80\%\\
\hline \first
Univ. of California-Santa Barbara& RU/VH& 149& Public& Far West& 21,026& 84\%\\
Univ. of California-Santa Cruz& RU/VH& 78& Public& Far West& 15,036& 88\%\\
Univ. of Chicago& RU/VH& 247& Private& Great Lakes& 13,870& 33\%\\
Univ. of Illinois-Chicago& RU/VH& 292& Public& Great Lakes& 24,865& 58\%\\
Univ. of Illinois-Urbana& RU/VH& 494& Public& Great Lakes& 40,687& 71\%\\
Univ. of Maryland-College Park& RU/VH& 322& Public& Mid East& 34,933& 68\%\\
Univ. of Michigan-Ann Arbor& RU/VH& 780& Public& Great Lakes& 39,533& 61\%\\
Univ. of Minnesota-Twin Cities& RU/VH& 509& Public& Plains& 50,954& 56\%\\
Univ. of N.\ Carolina-Chapel Hill& RU/VH& 391& Public& Southeast& 26,878& 59\%\\
Univ. of Notre Dame& RU/VH& 60& Private& Great Lakes& 11,479& 72\%\\
\hline \first
Univ. of Oregon& RU/H& 45& Public& Far West& 20,296& 76\%\\
Univ. of Pennsylvania& RU/VH& 565& Private& Mid East& 23,305& 46\%\\
Univ. of Southern California& RU/VH& 414& Private& Far West& 32,160& 50\%\\
Univ. of Texas-Austin& RU/VH& 344& Public& Southwest& 50,377& 70\%\\
Univ. of Utah& RU/VH& 224& Public& Rocky Mtns.& 28,933& 62\%\\
Univ. of Virginia-Charlottesville& RU/VH& 206& Public& Southeast& 23,341& 58\%\\
Univ. of Washington-Seattle& RU/VH& 685& Public& Far West& 39,199& 64\%\\
Univ. of Wisconsin-Madison& RU/VH& 717& Public& Great Lakes& 40,455& 68\%\\
Washington Univ.-St. Louis& RU/VH& 474& Private& Plains& 13,210& 49\%\\
Yale Univ.& RU/VH& 388& Private& New England& 11,441& 46\%\\
\hline \first
Aggregate && \multicolumn {1} {r} {\hspace*{-2em}17,381}&&& \multicolumn {1} {r} {\hspace*{-2em}1,227,269}& 63\%\\
\end {tabular}%
\end {center}%
\end {table}

\begin {table} 
\caption {Personnel of the institutions and rating of the mathematics departments. ``STEM'' is science, technology, engineering, and mathematics. ``NL'' are full-time non-ladder academic staff, ``-R'' research only. ``GS'' are graduate students including part-time  assistants. ``Rtng'' is the quality rating of mathematics research-doctorate programs \citep [appendix table H-4] {Goldberger1995}. ``Grp'' is the peer group of mathematics departments \citep {AMSsurvey} with group I consisting of those rating 3--5.00. Data in columns 2--5 are from the \citet {Carnegie2010}. Data in columns 6--8 were gathered in Spring 2010.}
\label {tab:personnel}
\begin {center}
\medskip
\scriptsize
\begin {tabular}
{l r r r r r r r r r}
\textsc {institution}& \multicolumn {2} {c} {Doctorates} &\multicolumn {2} {| c |} {All University}& \multicolumn {5} {c |} {Mathematics Department}\\ 
\cline {2-3} \cline {6-10}
\multicolumn {3} {r |} {\first STEM and/or Social Science Doctorates}&
\multicolumn {2} {c |} {Personnel}&
\multicolumn {3} {c |} {Personnel}&
\multicolumn {2} {c |} {Quality}\\
\cline {4-10}
\first & 
\multicolumn {1} {r |} {\hspace*{-108em}All Doctorates Annual Total}&
\multicolumn {1} {r |} {\footnotesize \%}&
\multicolumn {1} {c |} {NL-R}&
\multicolumn {2} {c |} {Ladder\hspace*{-0.5em}}&
NL&
\multicolumn {1} {r |} {GS}&
\multicolumn {1} {c} {Rtng.}&
\multicolumn {1} {r |} {\hspace*{-2em}Grp.\hspace*{-.25em}}\\
\hline \first
Auburn Univ.& 161& 31\%& 53& 1,088&&&& 2.31& II\\
Boston Univ.& 267& 57\%& 88& 1,314& 34& 41& 51& 3.03& I\\
Brandeis Univ.& 82& 52\%& 107& 347& 11& 6& 32& 3.64& I\\
Brown Univ.& 147& 69\%& 239& 778& 20& 13& 30& 3.73& I\\
California Inst. of Tech.& 166& 100\%& 527& 368& 15& 19& 32& 4.19& I\\
Carnegie Mellon Univ.& 195& 79\%& 248& 1,172& 38& 13& 42& 3.41& I\\
Columbia Univ.-New York City& 495& 52\%& 444& 3,221& 24& 42& 69& 4.23& I\\
Cornell Univ.& 412& 72\%& 749& 1,734& 45& 28& 68& 4.05& I\\
CUNY-Graduate School& 298& 51\%& 77& 133& 73& 0& 49& 3.65& I\\
Duke Univ.& 259& 68\%& 732& 2,902& 26& 34& 43& 3.53& I\\
\hline \first
Georgia Inst. of Tech.-Atlanta& 311& 96\%& 72& 856& 55& 29& 102& 3.19& I\\
Harvard Univ.& 572& 54\%& 3,950& 3,223& 21& 72& 58& 4.90& I\\
Indiana Univ.-Bloomington& 375& 37\%& 231& 1,486& 43& 21& 119& 3.53& I\\
Johns Hopkins Univ.& 362& 57\%& 1,421& 2,963& 16& 14& 32& 3.04& I\\
Massachusetts Inst. of Tech.& 467& 87\%& 1,026& 1,056& 51& 92& 90& 4.92& I\\
Michigan State Univ.& 430& 47\%& 393& 2,322& 60& 40& 113& 3.05& I\\
New York Univ.& 407& 40\%& 387& 3,083& 60& 109& 117& 4.49& I\\
Northwestern Univ.& 367& 57\%& 188& 1,996& 27& 23& 47& 3.71& I\\
Ohio State Univ.-Columbus& 560& 51\%& 414& 2,822& 64& 37& 197& 3.66& I\\
Penn.\ State Univ.-State College& 539& 60\%& 313& 2,817& 56& 41& 90& 3.50& I\\
\hline \first
Princeton Univ.& 276& 74\%& 381& 797& 36& 30& 69& \textbf {4.94}& I\\
Purdue Univ.-West Lafayette& 446& 59\%& 413& 2,060& 77& 66& 159& 3.82& I\\
Rensselaer Polytechnic Inst.& 156& 88\%& 83& 413& 27& 13& 60& 3.02& I\\
Rice Univ.& 126& 82\%& 147& 536& 15& 9& 30& 3.49& I\\
Rutgers Univ.-New Brunswick& 382& 62\%& 168& 2,158& 73& 77& 124& 3.96& I\\
Stanford Univ.& 625& 72\%& 1,503& 1,685& 28& 49& 74& 4.68& I\\
SUNY-Stony Brook& 285& 57\%& 229& 1,287& 32& 28& 71& 3.94& I\\
Univ. of California-Berkeley& 775& 62\%& 982& 1,495& 69& 34& 193& \textbf {4.94}& I\\
Univ. of California-Los Angeles& 666& 61\%& 1,481& 2,676& 64& 56& 200& 4.14& I\\
Univ. of California-San Diego& 327& 83\%& 1,246& 1,287& 58& 91& 131& 4.02& I\\
\hline \first
Univ. of California-Santa Barbara& 253& 69\%& 237& 907& 33& 25& 69& 3.04& I\\
Univ. of California-Santa Cruz& 107& 77\%& 113& 518& & & & 2.92& II\\
Univ. of Chicago& 331& 61\%& 546& 2,261& 37& 58& 87& 4.69& I\\
Univ. of Illinois-Chicago& 233& 44\%& 292& 1,925& 84& 3& 143& 3.58& I\\
Univ. of Illinois-Urbana& 574& 60\%& 279& 2,316& 70& 96& 181& 3.93& I\\
Univ. of Maryland-College Park& 482& 61\%& 244& 2,887& 59& 52& 216& 3.97& I\\
Univ. of Michigan-Ann Arbor& 660& 68\%& 909& 3,376& 69& 59& 143& 4.23& I\\
Univ. of Minnesota-Twin Cities& 592& 50\%& 935& 2,528& 68& 30& 132& 4.08& I\\
Univ. of N.\ Carolina-Chapel Hill& 439& 57\%& 660& 1,659& 37& 11& 66& 3.24& I\\
Univ. of Notre Dame& 149& 66\%& 150& 783& 44& 22& 67& 3.11& I\\
\hline \first
Univ. of Oregon& 164& 44\%& 104& 936& 28& 15& 55& 3.06& I\\
Univ. of Pennsylvania& 413& 52\%& 1,157& 1,671& 30& 29& 64& 3.97& I\\
Univ. of Southern California& 573& 35\%& 594& 2,001& 36& 17& 60& 3.23& I\\
Univ. of Texas-Austin& 702& 57\%& 304& 2,435& 61& 48& 134& 3.85& I\\
Univ. of Utah& 216& 61\%& 288& 2,238& 44& 52& 83& 3.52& I\\
Univ. of Virginia-Charlottesville& 358& 49\%& 425& 1,985& 26& 13& 50& 3.18& I\\
Univ. of Washington-Seattle& 503& 63\%& 1,123& 3,761& 62& 20& 97& 3.39& I\\
Univ. of Wisconsin-Madison& 628& 59\%& 842& 4,359& 61& 71& 131& 4.10& I\\
Washington Univ.-St. Louis& 241& 56\%& 778& 1,849& 25& 4& 33& 3.42& I\\
Yale Univ.& 332& 57\%& 1,135& 3,147& 16& 27& 35& 4.55& I\\
\hline \first
Aggregate &
\multicolumn {1} {r} {\hspace*{-2em}18,886}&
60\%&
\multicolumn {1} {r} {\hspace*{-2em}29,407}&
\multicolumn {1} {r} {\hspace*{-2em}93,617}&
\multicolumn {1} {r} {\hspace*{-2em}2,108}&
\multicolumn {1} {r} {\hspace*{-2em}1,779}&
\multicolumn {1} {r} {\hspace*{-2em}4,338}
\\
\end {tabular}%
\end {center}%
\end {table}

\newcommand {\twodept} [1] {\rotatebox {90} {\parbox {2em} {\setlength {\baselineskip} {0pt}\textsc {#1}}}}
\newcommand {\dept} [1] {\rotatebox {90} {\mbox {\textsc {#1}}}}
\newcommand {\squeez} [1] {\hspace*{-1em}#1\hspace*{-1em}}

\begin {table} 
\caption {University units (chiefly, departments) for which data are reported in Figures\ \ref {fig:productivityBYauthor} and \ref {fig:productivityBYdept}: ``\textsc {appl math}'' is applied mathematics, ``\textsc {bios}'' is biostatistics, ``\textsc {bus}'' is school of business, ``\textsc {cs(e)}'' is computer science and engineering, ``\textsc {e(c)e}'' is electrical and computer engineering, ``\textsc {eco}'' is economics, ``\textsc {epi}'' is epidemiology, ``\textsc {ie}'' is industrial engineering, ``\textsc {math}'' is mathematics, ``\textsc {me}'' is mechanical engineering, ``\textsc {or}'' is operations research, ``\textsc {phy}'' is physics, ``\textsc {stat}'' is statistics. The absence of an entry in this table indicates the university in question does not have the corresponding department.}
\label {tab:departments}
\begin {center}
\medskip
\scriptsize
\setlength {\tabcolsep} {0.5em}
\begin {tabular}
{l | c c | c c c c c c c c c l |}
\first \textsc {institution}
& \dept {appl math}
& \dept {math}
& \dept {bios, epi}
& \dept {bus}
& \dept {cs(e)}
& \dept {e(c)e}
& \dept {eco}
& \dept {me}
& \twodept {or,\\ ie}
& \dept {phy}
& \dept {stat}
& \multicolumn {1} {c |} {\textsc {other units}}
\\
\hline \first
Boston Univ.&&\chk & 2& \chk& \chk& \chk& \chk& \chk&& \chk &&\\
Brandeis Univ.&&\chk &&& \chk&& \chk& \chk&& \chk&&\\
Brown Univ.&\chk &\chk &&& \chk&& \chk&&& \chk&& Sch.\ of Eng.\\
California Inst. of Tech.&\chk &\chk &&& \chk& \chk&& \chk&& \chk&& Appl.\ Phy.\ Dept.\\
Carnegie Mellon Univ.&&\chk && \chk& \chk& \chk& \chk& \chk&& \chk& \chk& Soc.\ \& Dec.\ Sci.\ Dept.\\
Columbia Univ.-New York City&&\chk & 2& \chk& \chk& \chk& \chk& \chk& \chk& \chk& \chk& Appl.\ Math.\ \& Phy.\ Dept.\\
\hline \first
Cornell Univ. &&\chk &\chk&& \chk& \chk& \chk& \chk& \chk& \chk& \chk & Appl.\ Math.\ Ctr. \\
CUNY-Graduate School&&\chk &&\chk & 2& \chk& \chk& \chk&& \chk&& \\
Duke Univ.&&\chk && \chk& \chk& \chk& \chk& \chk&& \chk& \chk&\\
Georgia Inst. of Tech.-Atlanta&&\chk && \chk&& \chk& \chk& \chk& \chk& \chk&& Col.\ of Comp.\\
Harvard Univ.&&\chk & 2& \chk&&& \chk&&& \chk& \chk& Sch.\ of Eng.\\
Indiana Univ.-Bloomington&&\chk && \chk& \chk&& \chk&&& \chk& \chk&\\
\hline \first
Johns Hopkins Univ.&&\chk & 2& \chk& \chk& \chk& \chk& \chk&& \chk& \chk&\\
Massachusetts Inst. of Tech.&&\chk && \chk&& \chk& \chk& \chk&& \chk&&\\
Michigan State Univ.&&\chk & \chk& \chk&& \chk& \chk& \chk&& \chk& \chk&\\
New York Univ.&&\chk &&& \chk&& \chk&&& \chk&&\\
Northwestern Univ.&\chk &\chk && \chk&& \chk& \chk& \chk& \chk& \chk& \chk&\\
Ohio State Univ.-Columbus&&\chk && \chk& \chk& \chk& \chk& \chk& \chk& \chk& \chk&\\
\hline \first
Penn.\ State Univ.-State College&&\chk && \chk& \chk& \chk& \chk& \chk& \chk& \chk& \chk&\\
Princeton Univ.&& \chk&&& \chk& \chk& \chk& \chk& \chk& \chk&&Appl. \& Comp.\ Math.\ Pgm.\\
Purdue Univ.-West Lafayette&&\chk && \chk& \chk& \chk& \chk& \chk& \chk& \chk& \chk&\\
Rensselaer Polytechnic Inst.&& \chk&&& \chk& \chk& \chk& \chk& \chk& \chk&&\\
Rice Univ.& \chk& \chk&& \chk& \chk& \chk& \chk& \chk&& \chk&&\\
Rutgers Univ.-New Brunswick&& \chk&& \chk& \chk& \chk& \chk& \chk& 2& \chk& \chk&\\
\hline \first
Stanford Univ.&& \chk&& \chk& \chk& \chk& \chk& \chk& \chk& \chk& \chk& Appl.\ Phy.\ Dept.; ICME\\
SUNY-Stony Brook& \chk& \chk&& \chk& \chk& \chk& \chk& \chk& \chk& \chk&&\\
Univ. of California-Berkeley&& \chk& \chk& \chk&& \chk& \chk& \chk& \chk& \chk& \chk&\\
Univ. of California-Los Angeles&& \chk& 2& \chk& \chk& \chk& \chk& \chk&& \chk& \chk&\\
Univ. of California-San Diego&& \chk&& \chk& \chk& \chk& \chk& \chk&& \chk&& \\
Univ. of California-Santa Barbara&& \chk&&& \chk& \chk& \chk& \chk&& \chk&&\\
\hline \first
Univ. of Chicago&& \chk& \chk& \chk& \chk&&& \chk&& \chk& \chk &\\
Univ. of Illinois-Chicago&& \chk& \chk& \chk& \chk& \chk& \chk& \chk&& \chk&&\\
Univ. of Illinois-Urbana&& \chk&& \chk& \chk& \chk& \chk& \chk& \chk& \chk& \chk&\\
Univ. of Maryland-College Park&& \chk& \chk& \chk& \chk& \chk& \chk& \chk&& \chk&&\\
Univ. of Michigan-Ann Arbor&& \chk&& \chk&& \chk& \chk& \chk&& \chk& \chk&\\
Univ. of Minnesota-Twin Cities&& \chk& \chk& \chk& \chk& \chk& \chk& \chk&& \chk& \chk&\\
\hline \first
Univ. of N.\ Carolina-Chapel Hill&& \chk& 2& \chk& \chk&& \chk&&& \chk&& Stat.\ \& Op.\  Res.\ Dept.\\
Univ. of Notre Dame&& \chk&& \chk& \chk& \chk& \chk& \chk&& \chk&&\\
Univ. of Oregon&& \chk&& \chk& \chk&& \chk&&& \chk&& \\
Univ. of Pennsylvania&& \chk& \chk& \chk& \chk& \chk& \chk& \chk& \chk& \chk& \chk&\\
Univ. of Southern California&& \chk& \chk& \chk& \chk& \chk& \chk& \chk& \chk& \chk&&\\
Univ. of Texas-Austin&& \chk&& \chk& \chk& \chk& \chk& \chk&& \chk&& Stat.\ \& Sci.\ Comp.\ Div.\\
\hline \first
Univ. of Utah&& \chk& \chk& \chk& \chk& \chk& \chk& \chk&& \chk&&\\
Univ. of Virginia-Charlottesville&& \chk& \chk& \chk& \chk& \chk& \chk& \chk& \chk& \chk& \chk&\\
Univ. of Washington-Seattle& \chk& \chk& 2& \chk& \chk& \chk& \chk& \chk& \chk& \chk& \chk&\\
Univ. of Wisconsin-Madison&& \chk& \chk& \chk& \chk& \chk& \chk& \chk& \chk& \chk& \chk& Eng.\ Phy.\ Dept.\\
Washington Univ.-St. Louis&& \chk& \chk& \chk& \chk& \chk& \chk& \chk&& \chk&& Appl.\ Stat.\ Ctr.\\
Yale Univ.&& \chk& 2& \chk& \chk& \chk& \chk& \chk&& \chk& \chk&\\
\hline
\end {tabular}%
\end {center}%
\end {table}

\clearpage

\raggedright
%\bibliographystyle {plainnat}
%\bibliographystyle {spbasic}
%\bibliography {bibliography}

\end{document}